\documentclass{article}

\usepackage{amsmath}
\usepackage{amsfonts}
\usepackage{amssymb}
\usepackage{graphicx}  %
\usepackage{hyperref}
\usepackage{booktabs}
\usepackage{subfig}    %
\usepackage{caption}   %
\usepackage[utf8]{inputenc}
\usepackage{CJKutf8}

\usepackage{pdflscape}
\usepackage{multirow}
\usepackage{longtable}

\usepackage[capitalize]{cleveref}

\crefname{section}{Sec.}{Secs.}
\Crefname{section}{Section}{Sections}

\crefname{table}{Tab.}{Tabs.}      %
\Crefname{table}{Table}{Tables}    %

\crefname{figure}{Fig.}{Figs.}     %
\Crefname{figure}{Figure}{Figures} %

\usepackage[preprint]{neurips_2025}

\title{Understanding the Monty Hall Problem Through a Quantum Measurement Analogy}

\author{%
  Mo Li \\
  Independent Researcher \\
  \texttt{limo.research@gmail.com} \\
}

\begin{document}

\maketitle

\begin{abstract}
    The Monty Hall problem is a classic probability puzzle known for its counterintuitive solution, revealing fundamental discrepancies between mathematical reasoning and human intuition. To bridge this gap, we introduce a novel explanatory framework inspired by quantum measurement theory. Specifically, we conceptualize the host’s actions—opening doors to reveal non-prizes—as analogous to quantum measurements that cause asymmetric collapses of the probability distribution. This quantum-inspired interpretation not only clarifies why the intuitive misunderstanding arises but also provides generalized formulas consistent with standard Bayesian results. We further validate our analytical approach using Monte Carlo simulations across various problem settings, demonstrating precise agreement between theoretical predictions and empirical outcomes. Our quantum analogy thus offers a powerful pedagogical tool, enhancing intuitive understanding of conditional probability phenomena through the lens of probability redistribution and quantum-like measurement operations.
    
\end{abstract}

\section{Introduction}

The Monty Hall problem is a probability puzzle inspired by the American television game show \textit{Let's Make a Deal}. 
The standard version presents a scenario where a contestant faces three closed doors: behind one door is a car (the prize), and behind the other two are goats (non-prizes). After the contestant selects a door, the host—who knows the locations of all items—deliberately opens one of the remaining doors to reveal a goat, then offers the contestant the opportunity to switch their selection to the other unopened door. The counterintuitive result is that switching doors increases the winning probability from $1/3$ to $2/3$.

Despite its simple formulation, this problem has generated significant controversy, even among mathematicians and statisticians \citep{vos1990let}. The problem was first posed by \citet{selvin1975problem} and later popularized by Marilyn vos Savant's column in Parade magazine, where even many PhD mathematicians initially rejected the correct solution \citep{morgan1991monty}. The discrepancy between mathematical analysis and human intuition remains a fascinating subject in probability education. While conventional explanations rely on Bayesian analysis \citep{bayes1763problem} or frequentist arguments, we propose an alternative framework drawing inspiration from quantum measurement theory.

Our approach treats the set of doors as a quantum-like system where probabilities redistribute after measurement operations (the host opening a door). Rather than a physical quantum system, we use this as an interpretive device—a mental model that provides intuition about how information and probabilities evolve through the problem. This perspective is particularly illuminating when extending the problem to variants with many doors and prizes, where the informational asymmetry between contestant and host becomes more pronounced.

\section{The Standard Monty Hall Problem: A Quantum Perspective}

\subsection{Classical Formulation}

In the standard formulation, we have three doors with one prize. The contestant first selects a door, after which the host (who knows the prize location) opens one of the remaining doors to reveal a non-prize. The contestant then decides whether to stay with their original choice or switch to the remaining unopened door.

The classical analysis shows that:
\begin{itemize}
    \item If the contestant stays, their probability of winning is $1/3$ (their initial probability).
    \item If the contestant switches, their probability of winning is $2/3$.
\end{itemize}

This result follows from Bayes' theorem, but many find it unintuitive. As \citet{gill2011monty} argues, the key insight is not about probability per se, but about the precise mathematical modeling of the host's behavior.

\subsection{Quantum Measurement Analogy}

We propose viewing this problem through the lens of quantum measurement theory \citep{whitaker2020quantum}, where:

\begin{enumerate}
    \item The initial state is a superposition of all possible prize locations, with appropriate probabilities.
    \item The contestant's initial choice divides the system into ``selected door‘’ and ``other doors‘’ subspaces.
    \item The host's door-opening action represents a measurement that causes partial collapse of the probability wave function.
    \item Crucially, this measurement affects the ``other doors‘’ subspace differently than the ``selected door‘’ subspace.
\end{enumerate}

In the language of quantum mechanics, the host's action performs a non-uniform projection on the system's state, causing asymmetric probability redistribution. 

This is analogous to wavefunction collapse, where measurement of one part of an entangled system affects the probability distribution of the unmeasured parts.

Several researchers have explored formal quantum versions of the Monty Hall problem \citep{flitney2002quantum, dariano2002quantum, kurzyk2016quantum}, though our approach differs in that we use quantum concepts as an interpretive device for the classical problem rather than proposing a genuinely quantum-mechanical variant.

\section{Mathematical Analysis of Probability Transfer}

\subsection{General Formulation}

Let us develop a general framework for the Monty Hall problem with:

\begin{itemize}
    \item $N$ total doors
    \item $m$ doors with prizes
    \item The host opens $k$ doors
    \item Among the $k$ opened doors by the host, $r$ contain prizes
\end{itemize}

We denote:

\begin{itemize}
    \item $S_1$: The event that the contestant's initially selected door contains a prize
    \item $S_0$: The event that the contestant's initially selected door does not contain a prize
    \item $E$: The event that the host opens $k$ specific doors, revealing exactly $r$ prizes
\end{itemize}

\subsection{Bayesian Analysis with Different Host Strategies}

\subsubsection{Host with Perfect Information (Deliberate Selection)}

When the host knows the prize locations and deliberately selects which doors to open:

\begin{align}
P(S_1) &= \frac{m}{N} \\
P(S_0) &= \frac{N-m}{N} = 1 - P(S_1)
\end{align}

If the host always can and does reveal exactly $r$ prizes among the $k$ opened doors (event $E$), then:

\begin{align}
P(E|S_1) &= 1 \\
P(E|S_0) &= 1
\end{align}

Using Bayes' theorem, the posterior probability that the contestant's door contains a prize is:

\begin{equation} \label{eq:informed_host_stay}
\begin{split}
P(S_1|E) &= \frac{P(E|S_1)P(S_1)}{P(E|S_1)P(S_1) + P(E|S_0)P(S_0)} \\
&= \frac{1 \cdot \frac{m}{N}}{1 \cdot \frac{m}{N} + 1 \cdot (1-\frac{m}{N})} \\
&= \frac{m}{N}
\end{split}
\end{equation}

This shows that with a fully-informed host who deterministically opens doors, the probability of the contestant's original door containing a prize remains unchanged at $\frac{m}{N}$.
If the contestant switches to one of the remaining unopened doors, what is their probability of winning? When the contestant's original door contains a prize ($S_1$), there are $m-1-r$ prizes remaining among the $N-k-1$ unopened doors (since one prize is in their original door). When the contestant's original door does not contain a prize ($S_0$), there are $m-r$ prizes among the $N-k-1$ unopened doors. Therefore:

\begin{equation} \label{eq:informed_host_switch}
\begin{split}
P(\text{win by switching}) &= P(S_1)P(\text{win by switching}|S_1) + P(S_0)P(\text{win by switching}|S_0) \\
&= \frac{m}{N}\frac{m-1-r}{N-k-1} + (1-\frac{m}{N})\frac{m-r}{N-k-1} \\
&= \frac{m}{N}\frac{m-1-r}{N-k-1} + \frac{N-m}{N}\frac{m-r}{N-k-1} \\
&= \frac{1}{N(N-k-1)}[m(m-1-r) + (N-m)(m-r)] \\
&= \frac{1}{N(N-k-1)}[m^2-m-mr + Nm-rN-m^2+mr] \\
&= \frac{1}{N(N-k-1)}[Nm-m-rN] \\
&= \frac{1}{N(N-k-1)}[m(N-1)-rN] \\
&= \frac{m(N-1)-rN}{N(N-k-1)} \\
&= \frac{(m-r)-\frac{m}{N}}{N-k-1} \\
&= \frac{\frac{m}{N} \cdot (N-1) - r}{N-k-1}
\end{split}
\end{equation}

\subsubsection{Host with No Information (Random Selection)}

Now consider a host who randomly selects $k$ doors to open (excluding the contestant's door) and happens to reveal exactly $r$ prizes.

For $S_1$ (contestant's door has a prize), the probability of event $E$ is:
\begin{align}
P(E|S_1) &= \frac{\binom{m-1}{r}\binom{N-m}{k-r}}{\binom{N-1}{k}}
\end{align}

For $S_0$ (contestant's door has no prize), the probability of event $E$ is:
\begin{align}
P(E|S_0) &= \frac{\binom{m}{r}\binom{N-m-1}{k-r}}{\binom{N-1}{k}}
\end{align}

The posterior probability is:
\begin{align}
P(S_1|E) &= \frac{P(E|S_1)P(S_1)}{P(E|S_1)P(S_1) + P(E|S_0)P(S_0)} \\
&= \frac{\frac{m}{N}\frac{\binom{m-1}{r}\binom{N-m}{k-r}}{\binom{N-1}{k}}}{\frac{m}{N}\frac{\binom{m-1}{r}\binom{N-m}{k-r}}{\binom{N-1}{k}} + \frac{N-m}{N}\frac{\binom{m}{r}\binom{N-m-1}{k-r}}{\binom{N-1}{k}}} \\
&= \frac{m\binom{m-1}{r}\binom{N-m}{k-r}}{m\binom{m-1}{r}\binom{N-m}{k-r} + (N-m)\binom{m}{r}\binom{N-m-1}{k-r}} \\
&= \frac{m}{\displaystyle m + (N-m)\frac{\binom{m}{r}}{\binom{m-1}{r}}\frac{\binom{N-m-1}{k-r}}{\binom{N-m}{k-r}}}
\end{align}

Using the identities $\binom{m}{r} = \frac{m}{m-r}\binom{m-1}{r}$ and $\binom{N-m-1}{k-r} = \frac{N-m-k+r}{N-m}\binom{N-m}{k-r}$, we get:

\begin{equation} \label{eq:random_host_stay}
\begin{split}
P(\text{win by staying}) &= P(S_1|E) \\
&= \frac{m}{m + (N-m)\frac{m}{m-r}\frac{N-m-k+r}{N-m}} \\
&= \frac{m}{m + \frac{m(N-m-k+r)}{m-r}} \\
&= \frac{1}{1 + \frac{N-m-k+r}{m-r}} \\
&= \frac{m-r}{N-k}
\end{split}
\end{equation}

This elegant result has an intuitive interpretation: after revealing $r$ prizes among $k$ opened doors, the probability that the contestant's door contains a prize is simply the ratio of remaining prizes $(m-r)$ to the total number of unopened doors $(N-k)$.

To calculate the probability of winning by switching, we need to consider the posterior probabilities after the host has revealed $r$ prizes among $k$ opened doors (event $E$).

\begin{equation} \label{eq:random_host_switch}
\begin{split}
P(\text{win by switching}|E) &= \sum_{i=0}^{1} P(\text{win by switching}|S_i,E) \cdot P(S_i|E) \\
& = P(\text{win by switching}|S_1,E) \cdot P(S_1|E) \\
&\quad + P(\text{win by switching}|S_0,E) \cdot P(S_0|E)
\end{split}
\end{equation}

To derive this formula rigorously, we start from the definition of conditional probability:
\begin{align}
P(\text{win by switching}|E) = \frac{P(\text{win by switching} \cap E)}{P(E)}
\end{align}

Since $\{S_0, S_1\}$ forms a partition of the sample space (mutually exclusive and collectively exhaustive), we can decompose the numerator using the law of total probability:
\begin{align}
P(\text{win by switching} \cap E) &= P(\text{win by switching} \cap E \cap S_1) + P(\text{win by switching} \cap E \cap S_0)
\end{align}

Using the multiplication rule of probability, each term on the right-hand side can be rewritten as:
\begin{align}
P(\text{win by switching} \cap E \cap S_i) &= P(\text{win by switching}|E,S_i) \cdot P(E \cap S_i) \\
&= P(\text{win by switching}|E,S_i) \cdot P(S_i|E) \cdot P(E)
\end{align}

Substituting back and dividing both sides by $P(E)$:
\begin{align}
P(\text{win by switching}|E) &= \frac{P(\text{win by switching}|E,S_1) \cdot P(S_1|E) \cdot P(E)}{P(E)} \\
&\quad + \frac{P(\text{win by switching}|E,S_0) \cdot P(S_0|E) \cdot P(E)}{P(E)} \\
&= P(\text{win by switching}|E,S_1) \cdot P(S_1|E) \\
&\quad + P(\text{win by switching}|E,S_0) \cdot P(S_0|E)
\end{align}

This confirms our formula is correct. The key insight is that $\{S_0, S_1\}$ forms a partition of the sample space, and this property is preserved when intersected with any event $E$. The sets $S_0 \cap E$ and $S_1 \cap E$ remain mutually exclusive and together comprise the entire event $E$.

We've already calculated the posterior probabilities in equation \eqref{eq:random_host_stay}:
\begin{itemize}
    \item $P(S_1|E) = \frac{m-r}{N-k}$ (posterior probability the original door contains a prize)
    \item $P(S_0|E) = 1 - P(S_1|E) = \frac{N-m-k+r}{N-k}$ (posterior probability the original door does not contain a prize)
\end{itemize}

For the conditional probabilities of winning by switching:

\begin{itemize}
    \item If the original door has a prize ($S_1$), after the host opens $k$ doors revealing $r$ prizes, among the $N-k-1$ remaining unopened doors, there are $m-1-r$ prizes. Thus, $P(\text{win by switching}|S_1,E) = \frac{m-1-r}{N-k-1}$
    
    \item If the original door has no prize ($S_0$), after the host opens $k$ doors revealing $r$ prizes, among the $N-k-1$ remaining unopened doors, there are $m-r$ prizes. Thus, $P(\text{win by switching}|S_0,E) = \frac{m-r}{N-k-1}$
\end{itemize}

Substituting these values:

\begin{equation}
    \label{eq:random_host_switch_bayes_noinformation}
\begin{split}
P(\text{win by switching}|E) &= \frac{m-r}{N-k} \cdot \frac{m-1-r}{N-k-1} + \frac{N-m-k+r}{N-k} \cdot \frac{m-r}{N-k-1} \\
&= \frac{m-r}{(N-k)(N-k-1)}[(m-1-r) + (N-m-k+r)] \\
&= \frac{m-r}{(N-k)(N-k-1)}[N-k-1] \\
&= \frac{m-r}{N-k}
\end{split}
\end{equation}

Therefore, with a random host, the probability of winning by switching equals the probability of winning by staying: $\frac{m-r}{N-k}$. This means switching offers no advantage over staying.

This result makes intuitive sense: when the host randomly opens doors (without knowledge of the prize locations), the information revealed creates no asymmetry in the probability distribution across the unopened doors. Each unopened door, including the contestant's original choice, has the same probability of containing a prize, since their state is ``not having been observed''.
This contrasts sharply with the informed host scenario, where switching is advantageous because the host's knowledge creates an asymmetry in the probability distribution.

\section{Quantum Probability Flow }

We can understand the Monty Hall problem through what we term ``probability flow‘’ or ``probability transfer.‘’ Initially, probabilities are distributed across all doors according to their prior probabilities. When the host opens doors, revealing certain information, the probabilities associated with those doors must be redistributed—analogous to how quantum wave functions collapse upon measurement \citep{whitaker2020quantum}.

\subsection{The Flow of Probability}

In the standard three-door problem:
\begin{itemize}
    \item Initially, the contestant's door has $1/3$ probability of containing the prize.
    \item The other two doors collectively have $2/3$ probability.
    \item When the host (who knows the prize location) opens one of those doors to reveal a goat, the $2/3$ probability originally spread across two doors becomes concentrated in the single remaining unopened door.
\end{itemize}

This redistribution of probability—from opened doors to unopened ones—is what we call ``probability transfer.‘’ The key insight is that the transfer happens differently depending on the host's knowledge and door-selection strategy.

\subsection{Host Knowledge as Quantum Measurement}

The host's knowledge functions analogously to a quantum measurement device that collapses certain possibilities while preserving others. When the host deliberately avoids opening prize doors:

\begin{itemize}
    \item If the prize is behind the contestant's door (probability $1/3$), the host randomly opens one of the two goat doors.
    \item If the prize is behind one of the other doors (probability $2/3$), the host must open the specific goat door, leaving the prize door closed.
\end{itemize}

This asymmetry in the host's behavior creates the $2/3$ advantage for switching. The host's knowledge-guided action has effectively ``measured‘’ the system in a way that preserves information differently across the possible states.

This view aligns with the formal quantum treatment by \citet{flitney2002quantum}, who modeled the problem using quantum game theory, showing how quantum measurement concepts can clarify the classical solution.

\subsection{Generalized N-Door, M-Prize Case}

For the generalized case with $N$ doors and $m$ prizes, where the host opens $k$ doors revealing $r$ prizes:

\subsubsection{Host with No Information (Random Selection)}

When the host has no special knowledge and randomly opens doors, the quantum measurement analogy works differently:

Initially, all doors still have equal probability $\frac{m}{N}$ of containing prizes.

When the host randomly opens $k$ doors and reveals $r$ prizes, we need to consider how this affects the probability distribution across all doors. In this scenario, the host is truly selecting doors at random without any knowledge of prize locations.

Let's analyze this situation from the perspective of probability symmetry:

\begin{itemize}
    \item Initially, all $N$ doors have equal probability $\frac{m}{N}$ of containing a prize
    \item The contestant selects one door (which remains unopened)
    \item The host randomly opens $k$ doors from the remaining $N-1$ doors, revealing $r$ prizes
\end{itemize}

The key insight is that in this random selection scenario, all unopened doors share the same state of ``not having been observed.‘’ This creates symmetry between:
\begin{itemize}
    \item The contestant's initially chosen door
    \item The remaining unopened doors ($N-k-1$ doors)
\end{itemize}

Since the host's random selection doesn't create any asymmetry between the contestant's door and other unopened doors, they must all have equal probability of containing a prize. With $m-r$ prizes remaining behind $N-k$ unopened doors (including the contestant's door), each unopened door has probability:

\begin{align}
P(\text{prize in any unopened door}) = \frac{m-r}{N-k}
\end{align}

This means that switching doors offers no advantage. The probability of winning by staying with the original door equals the probability of winning by switching to any other unopened door.

However, it's important to note that the host's action of opening $k$ doors and revealing $r$ prizes does cause a redistribution of probabilities. Initially, the probability of winning by staying with the contestant's chosen door was $\frac{m}{N}$. After the host randomly opens doors, this probability actually changes. This change occurs because the host's measurement action breaks the uncertainty of the system. Once measurement occurs, it collapses the uncertainty and leads to a redistribution of probabilities across all unopened doors to $\frac{m-r}{N-k}$.

From a quantum measurement perspective, we can view this as all unopened doors sharing the same ``quantum state‘’ after measurement, with no door having privileged information compared to others. The measurement action only affects the doors that were actually opened, leaving all unopened doors in an equivalent probabilistic state.

\subsubsection{Host with Perfect Information}

When the host has perfect information about prize locations and deliberately chooses which doors to open, we can analyze the probability transfer using our quantum measurement framework:

Initially, all doors have equal probability of containing prizes, with each door having probability $\frac{m}{N}$ of containing a prize.

The key insight is that the host's deliberate selection process creates an asymmetry between the contestant's door and the other doors. The host effectively performs a measurement on all $N-1$ doors (excluding the contestant's choice), but only reveals $k$ of them. This measurement operation fundamentally alters the probability distribution.

We can conceptualize this process as follows:

\begin{itemize}
    \item The host observes all $N-1$ unchosen doors (introducing certainty through measurement)
    \item The host selectively reveals $k$ doors containing exactly $r$ prizes
    \item For the remaining unopened doors, a form of ``selective forgetting‘’ occurs - while the host knows their contents, from the contestant's perspective these doors return to a state of quantum-like superposition
\end{itemize}

Here, we can interpret the probability transfer process as a kind of conservation of expectation. 

That is, the expected number of prizes behind all doors (except the contestant's) is conserved throughout the host's measurement and revealing process. 

The contestant's door always retains its original probability of $\frac{m}{N}$ of containing a prize, since it is never measured or affected by the host's actions. 

For the other $N-1$ doors, before any doors are opened, the total expected number of prizes is $\frac{m}{N} \cdot (N-1)$. 

After the host deliberately opens $k$ doors and reveals $r$ prizes, the remaining $N-k-1$ unopened doors (excluding the contestant's) must contain the rest of the expected prizes, i.e., $\frac{m}{N} \cdot (N-1) - r$. 

This is a direct result of expectation conservation: the sum of the expected prizes behind the opened doors ($r$), the remaining unopened doors, and the contestant's door ($\frac{m}{N}$) must always add up to $m$.

Therefore, the probability of winning by switching to any one of the remaining unopened doors is simply the remaining expected number of prizes divided equally among those doors:
\begin{align}
P(\text{win by switching}) = \frac{\frac{m}{N} \cdot (N-1) - r}{N-k-1}
\end{align}
Note that this equation is the same as equation \ref{eq:informed_host_switch}.
In this sense, the ``probability flow'' from the opened doors to the remaining unopened doors can be understood as a redistribution of the expected value, with the total expectation conserved throughout the process. The contestant's door, never measured, keeps its original expectation, while the rest is shared among the other unopened doors.

The asymmetry arises because the host's measurement process affects only the unchosen doors, not the contestant's door. This is analogous to how in quantum mechanics, measuring one part of an entangled system affects only the measured components while leaving others unchanged.

This contrasts sharply with the random host scenario, where the host's actions create no asymmetry between the contestant's door and other unopened doors. In the informed host case, the ``measurement'' (opening doors) collapses the probability state asymmetrically, making switching advantageous.

\paragraph{How Many Doors Does the Host Actually Observe?}

In the case of a host with perfect information, a subtle but important question arises: \emph{How many doors does the host actually observe (i.e., measure) on stage?} This is not just a technicality, but is crucial for understanding the quantum analogy and the asymmetry in probability flow.

Suppose the host needs to select $k$ doors to open, such that exactly $r$ of them contain prizes. On the surface, it seems the host simply picks $k$ doors meeting this requirement. But in practice, to guarantee this selection, the host must know the state (prize or not) of all $N-1$ doors not chosen by the contestant. That is, the host must \emph{observe} (measure) all $N-1$ doors to find a subset of $k$ doors with exactly $r$ prizes.

This leads to a conceptual puzzle: if the host has measured all $N-1$ doors, why do we say that the probability associated with the $k$ opened doors is redistributed only among the remaining $N-k-1$ unopened doors? Why not among all $N-1$ doors, since all have been observed? What justifies the asymmetry?

The resolution is to recognize that, in our quantum-inspired framework, the host's knowledge is not retained after the selection. That is, after the host selects and opens $k$ doors (revealing $r$ prizes), we assume the host \emph{forgets} the detailed contents of the remaining $N-k-1$ unopened doors. This is analogous to a quantum measurement followed by selective ``forgetting": once the $k$ doors are opened, the information about the other doors is erased, and from the contestant's perspective, those doors return to a superposed (unmeasured) state.

This ``selective forgetting‘’ is a key postulate in our analysis: \textbf{if no observer retains information about a system, its state can be treated as a quantum-like superposition again.} Thus, only the $k$ opened doors are truly measured and collapsed; the rest revert to an unmeasured, symmetric state.

Therefore, the probability (or expected number of prizes) that was associated with the $k$ opened doors is redistributed \emph{only} among the $N-k-1$ remaining unopened doors. The contestant's door, never measured, keeps its original probability. This explains why, in the formula, the probability ``flows‘’ from the $k$ opened doors to the $N-k-1$ unopened ones, and not to all $N-1$ doors.

In summary, the host with perfect information must, in principle, observe all $N-1$ doors to select the $k$ to open, but after opening them, the information about the rest is ``forgotten,‘’ restoring symmetry and justifying the probability transfer in our equations.

\section{Concrete Examples}

Here we apply our derived formulas to specific instances of the Monty Hall problem, demonstrating how the quantum measurement framework provides consistent explanations across different scenarios.

\subsection{The Classic 3-Door, 1-Prize Problem}

For the classic Monty Hall scenario, we have $N=3$, $m=1$, $k=1$, $r=0$ (the host opens one door revealing no prizes). Let's apply our equations for both host types:

\subsubsection{Informed Host}

Using equation \ref{eq:informed_host_stay}, the probability of winning by staying with the original door is:
\begin{align}
P(\text{win by staying}) = \frac{m}{N} = \frac{1}{3}
\end{align}

For switching, we apply equation \ref{eq:informed_host_switch}:
\begin{align}
P(\text{win by switching}) &= \frac{\frac{m}{N} \cdot (N-1) - r}{N-k-1} \\
&= \frac{\frac{1}{3} \cdot (3-1) - 0}{3-1-1} \\
&= \frac{\frac{2}{3}}{1} = \frac{2}{3}
\end{align}

This confirms the well-known result that switching doors doubles the probability of winning when the host knows the prize location.

\subsubsection{Random Host}

For a host who randomly selects a door to open (excluding the contestant's door) and happens to reveal a non-prize, we apply equation \ref{eq:random_host_stay}:
\begin{align}
P(\text{win by staying}) &= \frac{m-r}{N-k} \\
&= \frac{1-0}{3-1} = \frac{1}{2}
\end{align}

The probability of winning by switching is, according to equation \ref{eq:random_host_switch_bayes_noinformation}:
\begin{align}
P(\text{win by switching}) &= \frac{m-r}{N-k} \\
&= \frac{1-0}{3-1} \\
&= \frac{1}{2}
\end{align}

With a random host, staying and switching yield equal probabilities of winning. This demonstrates why the host's knowledge is crucial to the puzzle's counterintuitive solution.

\subsection{A 100-Door, 37-Prize Problem}

Now we consider a more complex scenario with $N=100$ doors, $m=37$ prizes, where the host opens $k=3$ doors and reveals $r=2$ prizes.

\subsubsection{Informed Host}

For an informed host who deliberately reveals 2 prizes among the 3 opened doors, we use equation \ref{eq:informed_host_stay} to find:
\begin{align}
P(\text{win by staying}) = \frac{m}{N} = \frac{37}{100} = 0.37
\end{align}

And for switching, we apply equation \ref{eq:informed_host_switch}:
\begin{align}
P(\text{win by switching}) &= \frac{\frac{m}{N} \cdot (N-1) - r}{N-k-1} \\
&= \frac{\frac{37}{100} \cdot 99 - 2}{100-3-1} \\
&= \frac{36.63 - 2}{96} \\
&= \frac{34.63}{96} \approx 0.36073
\end{align}

Here, staying with the original door offers a slight advantage (37\% vs. 36.07\%).

\subsubsection{Random Host}

For a random host who happens to reveal 2 prizes among 3 randomly opened doors, we apply equation \ref{eq:random_host_stay}:
\begin{align}
P(\text{win by staying}) &= \frac{m-r}{N-k} \\
&= \frac{37-2}{100-3} \\
&= \frac{35}{97} \approx 0.36082
\end{align}

The probability of winning by switching (to any unopened door) is, by equation \ref{eq:random_host_switch_bayes_noinformation}:
\begin{align}
P(\text{win by switching}) &= \frac{m-r}{N-k} \\
&= \frac{37-2}{100-3} \\
&= \frac{35}{97} \approx 0.36082
\end{align}

However, since there are 96 doors to potentially switch to, the total probability of winning by randomly switching to any of the unopened doors is also $\frac{35}{97} \approx 0.3608$.

This example demonstrates that with many doors and prizes, the advantage of the informed host scenario diminishes compared to the classic 3-door case. It also shows how our quantum probability transfer framework consistently explains the probability distribution regardless of the problem's complexity.

\subsection{Numerical Simulation with Python}

To verify our analytical results, we used a Monte Carlo simulation implemented in Python to numerically study the generalized Monty Hall problem. We set the parameters as $N=100$ (number of doors), $m=37$ (number of prizes), $k=3$ (number of doors opened by the host), and $r=2$ (number of prizes revealed). The simulation considers both the \textbf{informed host} and \textbf{random host} scenarios, and records the empirical win rates for both ``stay'' and ``switch'' strategies.

\begin{table}[htbp]
    \centering
    \caption{Comparison of theoretical and simulated win rates for both host types and strategies. The differences are negligible, confirming the accuracy of our analytical results.}
    \begin{tabular}{lccc}
        \toprule
        & \textbf{Theoretical} & \textbf{Simulated} & \textbf{Difference} \\
        \midrule
        Informed host, Stay   & 0.37000 & 0.36982 & 0.00018 \\
        Informed host, Switch & 0.36073 & 0.36072 & 0.00001 \\
        Random host, Stay     & 0.36082 & 0.36097 & 0.00015 \\
        Random host, Switch   & 0.36082 & 0.36077 & 0.00005 \\
        \bottomrule
    \end{tabular}

\end{table}

\vspace{1em}

\begin{figure}[htbp]
    \centering
    
    \begin{minipage}{0.48\textwidth}
        \centering
        \includegraphics[width=\textwidth]{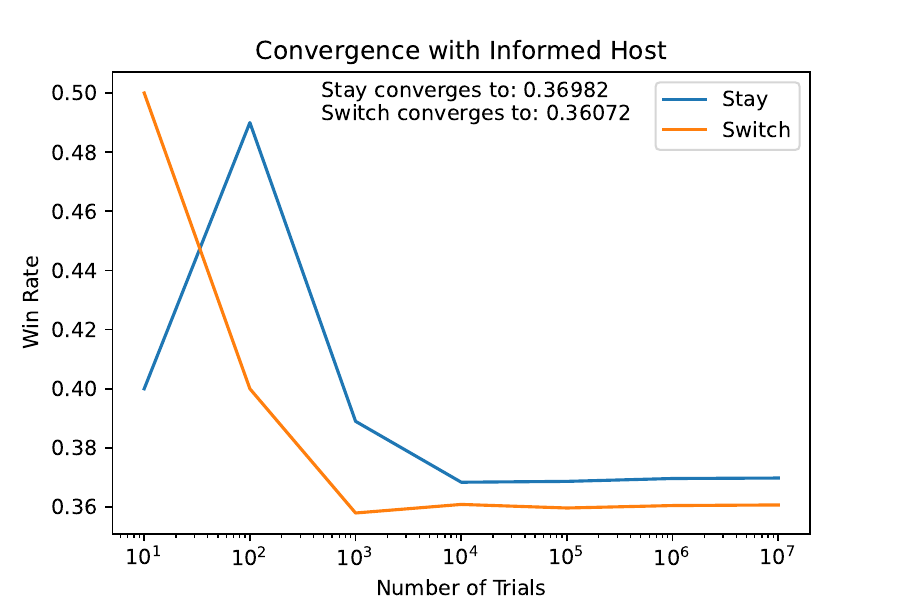}
    \end{minipage}
    \hfill
    \begin{minipage}{0.48\textwidth}
        \centering
        \includegraphics[width=\textwidth]{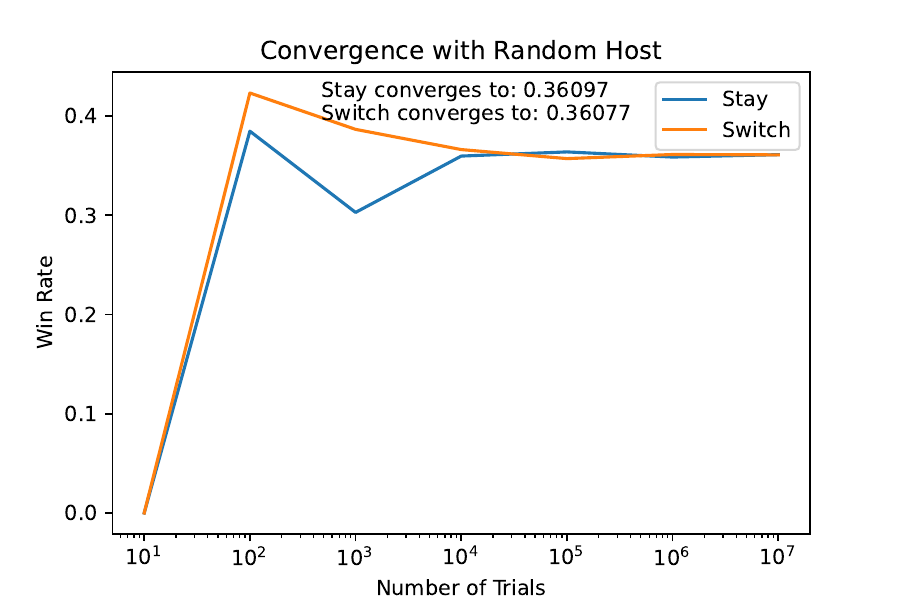}
    \end{minipage}
    \caption{Convergence of stay and switch win rates as the number of trials increases. Left: Informed host. Right: Random host.}
\end{figure}

These results further demonstrate that our quantum probability transfer framework and probability distribution analysis are consistent with empirical simulation.

\subsection{Quantum Measurement Interpretation}

These examples illustrate how quantum measurement concepts provide insight into conditional probability problems:

\begin{itemize}
    \item In the informed host scenario, the asymmetric measurement (deliberately revealing specific doors) creates a non-uniform redistribution of probability, affecting the optimal strategy.
    \item In the random host scenario, the measurement affects all unopened doors equally, leading to a uniform probability distribution among them.
    \item As the number of doors increases, the distinction between strategies becomes more subtle, showing how measurement effects scale with system size.
\end{itemize}

The quantum probability transfer framework thus provides a unified explanation for why the Monty Hall problem behaves differently under various host knowledge conditions, and how these effects manifest across different problem scales.

\section{Conclusion}

The quantum measurement perspective offers a powerful interpretive framework for the Monty Hall problem. By conceptualizing the host's actions as measurements that cause collapse and redistribution of probability, we gain intuition about the flow of information and probability in this and similar problems.
The key insight is that the host's knowledge and strategy determine how the measurement affects the probability distribution. When the host knows prize locations and deliberately avoids revealing prizes, the measurement preserves the probability of the contestant's door while concentrating the remaining probability into the other unopened doors. When the host has no special knowledge, the measurement affects all doors more uniformly.

This framework readily generalizes to variants of the problem with different numbers of doors and prizes, providing a unified understanding of conditional probability through the lens of quantum-inspired ``probability transfer''  mechanics. While not claiming any direct connection to physical quantum mechanics, this interpretive device bridges the gap between mathematical analysis and intuition in probability theory.
Our work builds upon and complements the formal quantum treatments by \citet{flitney2002quantum}, \citet{dariano2002quantum}, and \citet{kurzyk2016quantum}, but differs in focusing on the quantum analogy as a pedagogical and intuition-building device rather than as a genuinely quantum-mechanical alternative.

\begin{ack}
We thank Yuekai Li for discussions with the authors. We are also grateful for his helpful feedback on earlier drafts of this work.
\end{ack}

\bibliographystyle{plainnat}  %
\bibliography{references}

\appendix
\section{Related Work}

The Monty Hall problem has inspired numerous quantum extensions, leveraging quantum measurements and information concepts to reinterpret the classical paradox. Flitney and Abbott (2002) first introduced a quantum game-theoretic variant of the Monty Hall scenario, exploring quantum strategies and demonstrating how the host’s actions could be represented as quantum measurements, thus altering players’ optimal strategies and outcomes~\cite{FlitneyAbbott2002}. Building upon a similar quantum viewpoint, D’Ariano et al. (2002) proposed a continuous-state quantum Monty Hall game, examining optimal strategies when the host encodes his knowledge in classical versus quantum ways, notably illustrating advantages of quantum entanglement through a “quantum notepad” metaphor~\cite{dariano2002quantummontyhallproblem}.

Kurzyk and Glos (2016) provided another innovative angle by modeling the Monty Hall problem within quantum inferring acausal structures—quantum analogs of Bayesian networks—where conditional quantum operators replace classical conditional probabilities, thus highlighting the contrast between classical Bayesian updating and its quantum counterpart~\cite{Kurzyk_2016}. Ventura Olivella (2019) further introduced a quantum variant where the host lacks knowledge of the prize’s location and may accidentally reveal it, concluding that quantum strategies transform the originally biased game into a fair scenario, underscoring the role of measurement and information erasure in quantum probability redistribution~\cite{Olivella2019}.

Quezada et al. (2020) advanced this field practically by proposing an experimentally feasible quantum-optical setup to physically realize the Monty Hall game, examining the influence of entanglement, noise, and quantum channels on the game’s probability distribution, thereby offering empirical insight into quantum probability collapse and transfer mechanisms~\cite{Quezada_2020}. Lastly, Rajan and Visser (2019) connected the Monty Hall game with foundational quantum physics through the lens of the Pusey-Barrett-Rudolph (PBR) theorem, treating the classical puzzle as a thought experiment to explore epistemic versus ontic interpretations of quantum states, thus deepening the philosophical context for quantum measurement analogies~\cite{rajan2019quantumpbrtheoremmonty}.

In contrast to these works that pursue genuine quantum implementations or philosophical interpretations, our work introduces quantum measurement as an interpretative tool to elucidate the classical Monty Hall problem. We develop a unified quantum-inspired framework centered around probability collapse and selective information erasure, providing intuitive explanations for both classical and generalized variants involving arbitrary numbers of doors and prizes.

\section{Simulation Code}
\label{sec:simulation}

The core simulation code is available on GitHub Gist at the following link:

\begin{center}
\url{https://gist.github.com/Mor-Li/01788fb5006b8e4a57eee7324dd43409}
\end{center}

\noindent

\end{document}